\documentclass[11pt,openbib]{article}
\usepackage{latexsym}
\usepackage{amsmath}
\usepackage{amssymb}
\usepackage{amsfonts}
\numberwithin{equation}{section}
\newtheorem{thrm}{Theorem}[section]

\newtheorem{defn}{Definition}[section]
\newtheorem{lemm}{Lemma}[section]

\newtheorem{rmrk}{Remark}[section]

\newenvironment{prf}{{\em Proof:}}{\hfill{$\Box$}}
\newcommand{\A}{{\mathcal A}}
\newcommand{\alf}{\alpha}

\newcommand{\dt}[1]{\frac{d}{dt}{#1}}

\newcommand{\Dlt}{\Delta}

\newcommand{\field}[1]{\mathbb{#1}}

\newcommand{\Gmm}{\Gamma}

\newcommand{\gs}{\geqslant}

\newcommand{\half}{\frac{1}{2}}
\newcommand{\hf}{\half}

\newcommand{\lf}{\left}
\newcommand{\lra}[2]{\left\langle #1,\,#2\right\rangle}
\newcommand{\ls}{\leqslant}

\newcommand{\mpt}{\mapsto}

\newcommand{\nbl}{\nabla}

\newcommand{\ol}{\overline}
\newcommand{\omg}{\omega}
\newcommand{\Omg}{\Omega}
\newcommand{\prt}{\partial}
\newcommand{\pder}[2]{\frac{\partial#1}{\partial#2}}
\newcommand{\pdt}[1]{\frac{\partial#1}{\partial t}}
\newcommand{\pdz}[1]{\frac{\partial#1}{\partial z}}

\newcommand{\Di}{\prt_i}

\newcommand{\Px}{\prt_x}
\newcommand{\Py}{\prt_y}
\newcommand{\Pz}{\prt_z}

\newcommand{\R}{\field{R}}

\newcommand{\ra}{\rightarrow}

\newcommand{\rt}{\right}

\newcommand{\smgrp}{$\{S(t)\}_{t\ge 0}$}

\newcommand{\tld}{\widetilde}

\newcommand{\tht}{\theta}

\newcommand{\ups}{\upsilon}

\newcommand{\veps}{\varepsilon}

\newcommand{\vphi}{\varphi}
\newcommand{\vfi}{\varphi}

\newcommand{\zt}{\zeta}

\begin{document}
\pagestyle{myheadings}

\title{\bf On $H^2$ solutions and $z$-weak solutions of the 3D Primitive Equations}

\author{ Ning Ju $^1$}
\footnotetext[1]{
	Department of Mathematics, Oklahoma State University,
	401 Mathematical Sciences, Stillwater OK 74078, USA.
	Email: {\tt nju@okstate.edu},
	}
\date{July 21, 2015}

\maketitle

\begin{abstract}
 Global in time well-posedness of $H^2$ solutions and $z$-weak solutions of the
 3D Primitive equations in a bounded cylindrical domain is proved.
 More specifically, uniform in time boundedness and bounded absorbing sets are
 obtained for both $H^2$ solutions and $z$-weak solutions, as well as uniqueness
 of the $z$-weak solution for the 3D Primitive equations. The result for $H^2$
 solutions improves a recent one proved in \cite{j;t:2014}. The result for
 $z$-weak solution positively resolves the problem of global existence and
 uniqueness of $z$-weak solutions of the 3D primitive equations, which has been
 open since the work of \cite{tachim:2010}.
\\

\noindent
{\bf Keywords:} 3D viscous Primitive Equations, global existence, uniqueness,
 regularity.\\

{\bf MSC:} 35B40, 35B30, 35Q35, 35Q86.

\end{abstract}

%\tableofcontents

\markboth{3D Viscous Primitive Equations}
{\hspace{.2in} N. Ju \hspace{.5in} 3D Viscous Primitive Equations}

\indent
\baselineskip 0.58cm

\section{Introduction}
\label{s:intro}

 Given a bounded domain $D \subset \R^2$ with smooth boundary $\prt D$, we
 consider the following system of viscous Primitive Equations (PEs) of
 Geophysical Fluid Dynamics in the cylinder $\Omg =D\times(-h, 0)\subset\R^3$,
 where $h$ is a positive constant, see e.g. \cite{temma;ziane:04} and the
 references therein:

\noindent
 {\em Conservation of horizontal momentum:}
\begin{equation*}
  \label{e:v}
  \pdt{v} + (v\cdot\nbl)v + w \pdz{v} +\nbl p + fv^\bot + L_1 v = 0;
\end{equation*}
 {\em Hydrostatic balance:}
\begin{equation*}
  \label{e:hs}
  \Pz p + \tht =0 ;
\end{equation*}
 {\em Continuity equation:}
\begin{equation*}
  \label{e:cnt}
  \nbl\cdot v + \Pz w = 0 ;
\end{equation*}
 {\em Heat conduction:}
\begin{equation*}
  \label{e:t}
  \pdt{\tht} + v\cdot\nbl \tht + w \pdz{\tht} + L_2 \tht=Q.
\end{equation*}

 The unknowns in the above system of 3D viscous PEs are the fluid velocity field
 $(v, w)=(v_1, v_2, w)\in\R^3 $ with $v = (v_1, v_2)$ and $v^\bot=(-v_2, v_1)$
 being horizontal, the temperature $\tht$ and the pressure $p$. The Coriolis
 rotation frequency $f = f_0(\beta + y)$ in the $\beta$-plane approximation and
 the heat source $Q$ are given. For the issue concerned in this article, $Q$ is
 assumed to be independent of $t$. In the above equations and in this article,
 $\nbl$ and $\Dlt$ denote the horizontal gradient and Laplacian:
\[ \nbl := (\Px, \Py)  \equiv (\prt_1,\prt_2), \quad
   \Dlt := \Px^2+\Py^2 \equiv \sum_{i=1}^2\Di^2. \]
 We also use the following notation:
\[ \nbl_3 :=(\Px, \Py, \Pz)  \equiv (\prt_1,\prt_2, \prt_3). \]
 The viscosity and the heat diffusion operators $L_1$ and $L_2$ are given
 respectively as follows:
\begin{equation*}
 \label{e:L1.L2}
   L_i := - \nu_i\Dlt - \mu_i \frac{\prt^2}{\prt z^2}, \quad
	i =1, 2,
\end{equation*}
 where the positive constants $\nu_1, \mu_1$ are the horizontal and vertical
 viscosity coefficients and the positive constants $\nu_2, \mu_2$ are the
 horizontal and vertical heat diffusivity coefficients.

 The boundary of $\Omg$ is partitioned into three parts:
 $\partial\Omg = \Gmm_u \cup \Gmm_b \cup \Gmm_s$, where
\begin{align*}
\Gmm_u &:= \{(x, y, z) \in \overline{\Omg} : z = 0\},\\
\Gmm_b &:= \{(x, y, z) \in \overline{\Omg} : z = -h\},\\
\Gmm_s &:= \{(x, y, z) \in \overline{\Omg}:
 (x,y)\in \prt D\}.
\end{align*}

 Consider the following boundary conditions of the PEs as in \cite{cao;titi:05}
 and \cite{j:pe}:
\begin{equation*}
\begin{aligned}
  \mbox{on} \quad \Gmm_u \mbox{\,:}& \quad 
  \pdz{v}= h \tau, \quad w=0, \quad \pdz{\tht}=-\alf(\tht-\Theta),\\
  \mbox{on} \quad \Gmm_b \mbox{\,:}& \quad
  \pdz{v}=0,  \quad w=0, \quad \pdz{\tht}=0,\\
  \mbox{on} \quad \Gmm_s \mbox{\,:}& \quad 
  v\cdot n=0,  \quad \pder{v}{n}\times n= 0,
  \quad \pder{\tht}{n}=0,
\end{aligned}
\end{equation*}
 where $\tau(x, y)$ and $\Theta(x, y)$ are respectively the wind stress and
 typical temperature distribution on the surface of the ocean, $n$ is the
 normal vector of $\Gmm_s$ and $\alf$ is a non-negative constant.
 The above system of PEs will be solved with suitable initial conditions.

 We assume that $Q$, $\tau$ and $\Theta$ are independent of time. Notice that
 results similar to those to be presented here for the autonomous case can
 still be obtained for the non-autonomous case with proper modifications.
 For the autonomous case, assuming some natural compatibility conditions on
 $\tau$ and $\Theta$, one can further set $\tau=0$ and $\Theta=0$ without
 losing generality. See \cite{cao;titi:05} for a detailed discussion on this
 issue.

 Setting $\tau=0$, $\Theta=0$ and using the fact that
\begin{equation*}
 w(x, y, z, t) =  - \int_{-h}^z \nbl\cdot v(x, y, \xi, t)d\xi,
\end{equation*}
\begin{equation*}
 p(x, y, z, t) =  p_s(x, y, t) -\int^z_{-h} \tht(x, y, \xi, t)d\xi,
\end{equation*}
 one obtains the following equivalent formulation of the system of PEs:
\begin{equation}
\label{e:v.n}
\begin{split}
 \pdt{v} + L_1v 
 &+ (v\cdot\nbl)v
 -\left(\int_{-h}^z\nbl\cdot v(x,y,\xi,t)d\xi\right)\pdz{v}\\
 &+ \nbl p_s(x,y,t)-\int_{-h}^z \nbl\tht(x,y,\xi,t)d\xi
 + fv^\bot =0.
\end{split}
\end{equation}
\begin{equation}
\label{e:t.n}
 \pdt{\tht} + L_2\tht + v\cdot\nbl\tht 
 - \left(\int_{-h}^z\nbl\cdot v(x,y,\xi,t)d\xi\right)\pdz{\tht}
 =Q ;
\end{equation}
\begin{equation}
\label{e:bc.v}
  \pdz{v} \Big\vert_{z=0} = \pdz{v}\Big|_{z=-h} =0, \quad
 v\cdot n\big|_{\Gmm_s} =0, \quad \pder{v}{n}\times n\Big|_{\Gmm_s} =0,
\end{equation}
\begin{equation}
\label{e:bc.t}
 \lf(\pdz{\tht}+\alf \tht\rt)\Big|_{z=0}= \pdz{\tht}\Big|_{z=-h}=0, \quad
 \pder{\tht}{n}\Big|_{\Gmm_s} =0,
\end{equation}
\begin{equation}
\label{e:ic.n}
 v(x, y, z, 0) = v_0(x, y, z), \quad \tht(x, y, z, 0) = \tht_0(x, y, z).
\end{equation}

 We remark that the expressions of $w$ and $p$ via integrating the continuity
 equation and the hydrostatic balance equation were already used in
 \cite{lions;temam;wang:92.2} dealing with the Primitive Equations for large
 scale oceans. See also \cite{lions;temam;wang:92.1} for a similiar treatment
 of the Primitive Equations for atmosphere.

 Notice that the effect of the {\em salinity} is omitted in the above 3D
 viscous PEs for brevity of presentation. However, our results in this article
 are still valid when the effect of salinity is included. For simplicity of
 discussion, we only consider the case of $Q\in L^2$ being independent of time
 and set the right-hand side of \eqref{e:v.n} as zero. This is not techically
 essential. If the right-hand side of \eqref{e:v.n} is replaced by a non-zero
 time-independent given external force $R\in L^2(\Omega)$, the results of this
 paper are still valid. The case with time-dependent $Q$ and $R$ can be treated
 similarly with minor proper adjustments.

 To the best of our knowledge, the mathematical framework of the viscous
 primitive equations for the large scale ocean was first formulated in
 \cite{lions;temam;wang:92.2}; the notions of weak and strong solutions were
 defined and existence of weak solutions was proved. Uniqueness of weak
 solutions is still unresolved yet.
 Existence of strong solutions {\em local in time} and their uniqueness were
 obtained in \cite{g;m;r:01} and \cite{temma;ziane:04}. Existence of strong
 solutions {\em global in time} was proved independently in \cite{cao;titi:05}
 and \cite{kob:07}. See also \cite{kz:07} for dealing with some other boundary
 conditions. In \cite{j:pe}, existence of the global attractor for the strong
 solutions of the system is proved in the functional space of strong solutions.
 It is proved in \cite{j;t:2014} that the $H^2$ solutions are uniformly bounded
 with a bounded absorbing set in $H^2$ space under the assumption that $\alf=0$
 and $Q, Q_z\in L^2$ and thus the global attractor of the strong solutions has
 finite Hausdorff and fractal dimensions. In \cite{j:pedim}, it is further
 proved that the global attractor of the strong solutions has finite Hausdorff
 and fractal dimensions for any $\alf\gs0$ and $Q\in L^2$ via an approach
 different from that of \cite{j;t:2014}, without using uniform boundedness of
 $H^2$ solutions.

 This article focus on the study of global existence of $H^2$ solutions and
 $z$-weak solutions of the system of 3D viscous PEs and the uniqueness of
 $z$-weak solutions. Both of the problem are resolved in the sense to be
 discussed next.

 The global existence of the $H^2$ solutions was obtained in \cite{petcu:06}
 for the 3D viscous PEs with {\em periodic} boundary conditions,
 under the condition that $Q\in H^1$. However, this approach is not applicable
 to the case of non-periodic boundary conditions. A different and more involved
 analysis was presented recently in \cite{j;t:2014}, which proves the uniform
 boundedness and existence of a bounded absorbing set for the $H^2$ solutions
 of the 3D viscous PEs with the set of boundary conditions as given by
 \eqref{e:bc.v} and \eqref{e:bc.t} under even less demanding condition that
 $Q, Q_z\in L^2$. This analysis also applies to the case with
 {\em periodic} boundary conditions, thus improving the result of
 \cite{petcu:06}. However, the result of \cite{j;t:2014} requires the
 conditions that $Q_z\in L^2$ and that $\alf=0$. As the first main result of
 this article, we further eliminate these two extra conditions via a somewhat
 different approach. See Theorem~\ref{t:h2} in Section~\ref{s:h2}. The main
 idea relies on a new a priori esitmate for $(v_t,\tht_t)$ that was obtained
 recently in \cite{j:pedim}.

 The notion of a $z$-weak solution, i.e. a weak solution $(v,\tht)$ such that
\[ (v_z,\tht_z)\in L^\infty(0,T; H)\cap L^2(0,T; V), \]
 was introduced in \cite{petcu:07}, which proved global existence and
 uniqueness of $z$-weak solutions in a 2D domain with periodic boundary
 conditions.
 This problem was later studied in \cite{tachim:2010} for a 3D domain with
 non-periodic boundary conditions. Global existence and uniqueness of $z$-weak
 solutions was proved in \cite{tachim:2010} under the extra condition that
 $(v_0,\tht_0)\in L^6$.
 It seems to have been an open problem since the work of \cite{tachim:2010}
 that whether or not global existence and uniqueness of $z$-weak solutions
 are still valid for the 3D case. As our second main result, we resolve this
 problem positively.
 See Theorem~\ref{t:zw.e} in Section~\ref{s:zw.e} and Theorem~\ref{t:zw.u} in
 Section~\ref{s:zw.u}.

 The rest of this article is organized as follows:

 In Section~\ref{s:pre}, we give the notations, briefly review the background
 results and present the problems to be studied and recall some important facts
 crucial to later analysis.
 In Section~\ref{s:h2}, we prove our first main result, Theorem~\ref{t:h2}, on
 global existence and uniform boundedness of $H^2$ solutions for any $\alf\gs0$
 and $Q\in L^2$.
 In Section~\ref{s:zw.e}, we prove Theorem~\ref{t:zw.e} on global existence and
 uniform boundedness of $z$-weak solutions.
 In Section~\ref{s:zw.u}, we prove Theorem~\ref{t:zw.u} on uniqueness of
 $z$-weak solutions.

\section{Preliminaries}
\label{s:pre}

\noindent

 We recall that $D$ is a bounded smooth domain in $\R^2$ and
 $\Omg=D\times[0,-h]$, where $h$ is a positive constant. We denote by
 $L^p(\Omg)$ and $L^p(D)$ ($1\ls p <+\infty$) the classic $L^p$ spaces with
 the norms:
 \begin{equation*}
 \|\phi\|_p = \left\{ \begin{array}{ll}
   \left( \int_\Omg|\phi(x,y,z)|^pdxdydz\right)^\frac{1}{p},
	& \forall \phi\in L^p(\Omg); \\
  \left( \int_D|\phi(x,y)|^pdxdy \right)^\frac{1}{p},
	& \forall \phi\in L^p(D).
  \end{array} \right.
 \end{equation*}
 Denote by $H^m(\Omg)$ and $H^m(D)$ ($m\gs 1$) the classic Sobolev spaces for
 square-integrable functions with square-integrable derivatives up to order $m$.
 We do not distinguish the notations for vector and scalar function spaces,
 which are self-evident from the context. For simplicity, we may use $d\Omg$ to
 denote $dxdydz$ and $dD$ to denote $dxdy$, or we may simply omit them when
 there is no confusion. Using the H\"older inequality, it is easy to show that, 
 for $\vphi\in L^p(\Omg)$,
\begin{equation*}
%\label{e:lp.av}
 \|\ol{\vphi}\|_{L^p(\Omg)}  =  h^\frac{1}{p}\|\ol{\vphi}\|_{L^p(D)}
 \ls 
 \|\vphi\|_p, \quad \forall p\in[1, +\infty],
\end{equation*}
 where $\ol{\vphi}$ is defined as the vertical average of $\vphi$:
\[ \ol{\vphi}(x,y) = h^{-1}\int_{-h}^0 \vphi(x,y,z)dz .\]

 Define the function spaces $H$ and $V$ as follows:
\begin{equation*}
\begin{split}
 H := H_1\times H_2 := \{ v\in L^2(\Omg)^2\ | \
  \nbl\cdot\ol{v} =0,\quad \ol{v}\cdot n |_{\Gmm_s}=0 \}\times L^2(\Omg),\\
 V := V_1\times V_2 := \{ v\in H^1(\Omg)^2\ | \
  \nbl\cdot\ol{v} =0,\quad v\cdot n |_{\Gmm_s}=0 \} \times H^1(\Omg).
\end{split}
\end{equation*}
 Define the bilinear forms: $a_i: V_i\times V_i \ra \R$, $i=1,2$ as follows:
\begin{equation*}
\begin{aligned}
  a_1(v, u) &= \int_\Omg \left(\nu_1 \nbl v_1\cdot \nbl u_1 
		+\nu_1 \nbl v_2\cdot \nbl u_2
		+ \mu_1 v_z\cdot u_z \right)d\Omg;\\
  a_2(\tht, \eta)& = \int_\Omg \left(\nu_2 \nbl\tht \cdot \nbl\eta 
		+ \mu_2 \tht_z \eta_z \right)d\Omg
                + \alf\int_{\Gmm_u}\tht\eta dxdy.
\end{aligned}
\end{equation*}
 Let $V_i'$ ($i=1,2$) denote the dual space of $V_i$. We define the linear
 operators $A_i : V_i \mapsto V_i'$, $i=1,2$ as follows:
 \begin{equation*}
 \lra{A_1v}{u}=a_1(v, u), \quad \forall v, u \in V_1;
 \quad
 \lra{A_2\tht}{\eta}=a_2(\tht, \eta),
 \quad  \forall \tht, \eta \in V_2,
 \end{equation*}
 where $\lra{\cdot}{\cdot}$ is the corresponding scalar product between $V_i'$
 and $V_i$. We also use $\lra{\cdot}{\cdot}$ to denote the inner products in
 $H_1$ and $H_2$. Define:
 \[ D(A_i) = \{ \phi \in V_i, A_i\phi\in H_i \},\quad i=1,2. \]
 Since $A_i^{-1}$ is a self-adjoint compact operator in $H_i$, by the classic
 spectral theory, the power $A_i^s$ can be defined for any $s\in\R$.
 Then $D(A_i)' = D(A_i^{-1})$ is the dual space of $D(A_i)$ and
 $ V_i = D(A_i^\half)$, $ V_i' = D(A_i^{-\half})$. Moreover,
  \[ D(A_i) \subset V_i \subset H_i \subset V_i' \subset D(A_i)', \]
 where the embeddings above are all compact.
 Define the norm $\|\cdot\|_{V_i}$ by:
 \[ \|\cdot\|_{V_i}^2 = a_i(\cdot,\cdot)
	= \lra{A_i\cdot}{\cdot}
	= \lra{A_i^\half\cdot}{A_i^\half\cdot},
	\quad i=1,2. \]
 The Poincar\'e inequalities are valid. There is a constant $c>0 $,
 such that for any $\phi =(\phi_1, \phi_2)\in V_1$ and $\psi \in V_2$
\begin{equation*}
 c\|\phi\|_2 \ls \|\phi\|_{V_1}, \quad c\|\psi\|_2 \ls \|\psi\|_{V_2}.
\end{equation*}
 Therefore, there exist constants $c>0$ and $C>0$ such that for any
 $\phi =(\phi_1, \phi_2)\in V_1$ and $\psi \in V_2$,
\begin{equation*}
\label{e:nmeqv}
 c\|\phi\|_{V_1} \ls \|\phi\|_{H^1(\Omg)} \ls C \|\phi\|_{V_1},
\quad
 c\|\psi\|_{V_2} \ls \|\psi\|_{H^1(\Omg)} \ls C \|\psi\|_{V_2}.
\end{equation*}
 Notice that, in the above first inequality, we have written
 $\|\phi\|_{H^1(\Omg)}$ instead of $\|\phi\|_{H^1(\Omg)^2}$. We could also simply
 write $\|\phi\|_{H^1}$. Here and later on as well, we do not distinguish the
 notations for vector and scalar function spaces which are self-evident from
 the context. In this article, we use $c$ and $C$ to denote generic positive
 constants, the values of which may vary from one place to another.

 Recall the following definitions of weak and strong solutions:
\begin{defn}
\label{d:soln}
 Suppose $Q \in L^2(\Omg)$, $(v_0, \tht_0)\in H$ and $T>0$. The pair
 $(v, \tht)$ is called a {\em weak solution} of the 3D viscous PEs
 (\ref{e:v.n})-(\ref{e:ic.n}) on the time interval $[0, T]$ if it satisfies
 (\ref{e:v.n})-(\ref{e:t.n}) in the weak sense, and also
\begin{equation*}
 (v, \tht) \in C([0, T];H) \cap L^2(0, T; V), \quad
 \prt_t(v,\tht) \in L^1(0, T; V').
\end{equation*}
 If $(v_0, \tht_0)\in H$ and $(\prt_zv_0, \prt_z\tht_0)\in H$, a weak solution
 $(v, \tht)$ is called a {\em $z$-weak solution} of (\ref{e:v.n})-(\ref{e:ic.n})
 on the time interval $[0, T]$ if, in addition, it satisfies
\[ (v_z,\tht_z) \in L^\infty([0, T]; H) \cap L^2(0, T; V). \]
 Moreover, if $(v_0, \tht_0)\in V$, a weak solution $(v, \tht)$ is called a
 {\em strong solution} of (\ref{e:v.n})-(\ref{e:ic.n}) on the time interval
 $[0, T]$ if, in addition, it satisfies
\begin{equation*}
 (v,\tht) \in C([0, T]; V) \cap L^2(0, T; D(A_1)\times D(A_2)).
\end{equation*}

\end{defn}

 The following theorem on global existence and uniqueness for the strong
 solutions was proved in \cite{cao;titi:05}. See also a related result in
 \cite{kob:07}.

\begin{thrm}
\label{t:strong.g}
 Suppose $Q\in H^1(\Omg)$. Then, for every $(v_0,\tht_0)\in V$ and $T>0$, there
 exists a unique strong solution $(v, \tht)$ on $[0, T]$ to the system of 3D
 viscous PEs, which depends on the initial data continuously in $H$.
\end{thrm}

\begin{rmrk}
 It is easy to see from the proof of Theorem~\ref{t:strong.g} given in
 \cite{cao;titi:05} that the condition $Q\in H^1(\Omg)$ can be relaxed to
 $Q\in L^6(\Omg)$. Notice that there are gaps between Definition \ref{d:soln}
 and Theorem~\ref{t:strong.g} for the condition on $Q$, for the continuity
 of the strong solution with respect to time and for the continuous dependence
 of the strong solution with respect to initial data.
\end{rmrk}

 We now recall the following result proven in \cite{j:pe} for the existence
 of global attractor $\A$ for the strong solutions of the 3D viscous PEs
 (\ref{e:v.n})-(\ref{e:ic.n}).
\begin{thrm}
\label{t:attractor}
 Suppose that $Q\in L^2(\Omg)$ is independent of time. Then the solution
 operator \smgrp \ of the 3D viscous PEs (\ref{e:v.n})-(\ref{e:ic.n}):
 $S(t)(v_0, \tht_0) =(v(t), \tht(t))$ defines a semigroup in the space
 $V$ for $t\in \R_+$.  Moreover, the following statements are valid:
 \begin{enumerate}
 \item
 For any $(v_0, \tht_0)\in V$, $t\mpt S(t)(v_0,\tht_0)$ is continuous from
 $\R_+$ into $V$.
 \item
 For any $t>0$, $S(t)$ is a continuous and compact map in $V$.
 \item
 \smgrp \ possesses a global attractor $\A$ in the space $V$. The global
 attractor $\A$ is compact and connected in $V$ and it is the minimal bounded
 attractor in $V$ in the sense of the set inclusion relation; $\A$ attracts all
 bounded subsets of $V$ in the norm of $V$.
 \end{enumerate}
\end{thrm}

 We recall also the following important result proved in \cite{j:pedim}:
\begin{thrm}
\label{t:t}
 Suppose $Q \in L^2(\Omg)$ and $\alf\gs0$.

 For any $(v_0,\tht_0)\in V$ and $(\prt_tv(0),\prt_t\tht(0))\in H$, there
 exists a unique solution
 $(v,\tht)$ of \eqref{e:v.n}-\eqref{e:ic.n} such that
\begin{equation*}
 (\prt_tv,\prt_t\tht)\in L^\infty(0,\infty; H)\cap L^2(0,\infty;V).
\end{equation*}
 Moreover, there exists a bounded absorbing ball for $(\prt_tv,\prt_t\tht)$ in
 space $H$.
\end{thrm}

 The main goal of this article is to prove global in time uniform boundedness
 of the $H^2$ solutions and global existence and uniqueness of $z$-weak or the
 so-called vorticity solutions, we also prove existence of bounded absorbing
 set of $z$-weak solutions.

 To end this section, we introduce two technically critical lemmas to be used
 in our later analysis.
 First, we recall a lemma which will be useful for the {\em a priori} estimates
 in Section~\ref{s:h2}. See \cite{cao;titi:03}, and also \cite{j:pe}, for a
 proof.
\begin{lemm}
\label{l:ju.1}
 Suppose that $\nbl\ups, \vfi\in H^1(\Omg), \psi\in L^2(\Omg)$.
 Then, there exists a constant $C>0$ independent of $\ups, \vfi, \psi$ and $h$,
 such that
\begin{equation*}
\left|
 \lra{\left(\int_{-h}^z\nbl\cdot \ups(x,y,\xi)d\xi\right) \vfi}{\psi}
 \right|
 \ls C \|\nbl\ups\|_2^{\hf}\left\|\nbl \ups\right\|_{H^1}^{\hf}
	\|\vfi\|_{2}^{\hf}\|\vfi\|_{H^1}^{\hf} \|\psi\|_{2}.
\end{equation*}

\end{lemm}

 Next, we prove a different anisotropic estimate which will be crucial in the
 analysis of Section~\ref{s:zw.e} and Section~\ref{s:zw.u}.

\begin{lemm}
\label{l:ju.2}
 Suppose that $\Omg=D\times [-h,0]$ as given previously and $\phi$, $\psi$,
 $\vfi$, $\nbl \phi$, $\nbl \psi$, $\vfi_z\in L^2(\Omg)$. Then, there exists
 a positive constant $C$ which is independent of $\phi$, $\psi$ and $\vfi$,
 such that
\begin{equation*}
\begin{split}
 \int_\Omg |\phi\psi\vfi| \ls& C \|\phi\|_2^\hf(\|\phi\|_2+\|\nbl\phi\|_2)^\hf
 \|\psi\|_2^\hf(\|\psi\|_2+\|\nbl\psi\|_2)^\hf \\
 &\times \|\vfi\|_2^\hf(\|\vfi\|_2+\|\vfi_z\|_2)^\hf.
\end{split}
\end{equation*}

\end{lemm}

\begin{prf} By Cauchy-Schwartz inequality, we have
\begin{equation*}
\begin{split}
 \int_\Omg |\phi\psi\vfi|
&\ls \int_D \|\vfi\|_{L^\infty_z} \int_{-h}^0 |\phi\psi|\ dz dD\\
&\ls \lf(\int_D \|\vfi\|_{L^\infty_z}^2\rt)^\hf
     \lf[\int_D \lf(\int_{-h}^0 |\phi\psi|\ dz\rt)^2\rt]^\hf.
\end{split}
\end{equation*}
 By Agmon's inequality and Cauchy-Schwartz inequality, we have
\begin{equation*}
\begin{split}
 \int_D \|\vfi\|_{L^\infty_z}^2 
&\ls C \int_D \lf(\int_{-h}^0 |\vfi|^2\rt)^\hf
       \lf[\int_{-h}^0 \lf(|\vfi|^2 + |\vfi_z|^2\rt) \rt]^\hf \\
&\ls C\|\vfi\|_2 (\|\vfi\|_2+ \|\vfi_z\|_2).
\end{split}
\end{equation*}
 By Minkowski's inequality and Cauchy-Schwartz inequality, we have
\begin{equation*}
\begin{split}
\lf[\int_D \lf(\int_{-h}^0 |\phi\psi|\ dz\rt)^2\rt]^\hf
&\ls \int_{-h}^0 \lf(\int_D|\phi\psi|^2 \rt)^\hf \  dz \\
&\ls \int_{-h}^0 \lf(\int_D|\phi|^4\rt)^{\frac{1}{4}}
 \lf(\int_D|\psi|^4\rt)^{\frac{1}{4}} \ dz \\
&\ls C\int_{-h}^0 \|\phi\|_{L^2(D)}^\hf
 (\|\phi\|_{L^2(D)}+\|\nbl\phi\|_{L^2(D)})^\hf\ \\
 &\text{\hspace{.5in}}\times \|\psi\|_{L^2(D)}^\hf
 (\|\psi\|_{L^2(D)}+\|\nbl\psi\|_{L^2(D)})^\hf \  dz \\
&\ls C\|\phi\|_2^\hf(\|\phi\|_2+\|\nbl\phi\|_2)^\hf
      \|\psi\|_2^\hf(\|\psi\|_2+\|\nbl\psi\|_2)^\hf,
\end{split}
\end{equation*}
 where in the above second to the last step we have used
 Gagliardo-Nirenberg-Sobolev inequality. This finishes the proof.

\end{prf}

\section{The bounded absorbing set in $D(A_1)\times D(A_2)$}
\label{s:h2}

 The existence of solutions, global in time, in the space $H^2$ was first proved
 in \cite{j;t:2014} for the case with boundary conditions \eqref{e:bc.v} and
 \eqref{e:bc.t}. Uniform boundedness of the solutions in the space
 $D(A_1)\times D(A_2)$ and the existence of a bounded absorbing ball for the
 solutions in the space $D(A_1)\times D(A_2)$ were obtained in \cite{j;t:2014}
 under the condition that $\alf=0$ and that $Q, Q_z\in L^2$. What will be
 established in the following Theorem~\ref{t:h2} is that the above result is
 still valid for any $\alf\gs0$ and $Q\in L^2$, thus eliminating the extra
 conditions $\alf=0$ and $Q_z\in L^2$. The proof to be presented here will be
 different from the one given by \cite{j;t:2014}.

 Now, we state our first main result of this article in the following theorem:
\begin{thrm}
\label{t:h2}
 Suppose $Q\in L^2(\Omg)$. For any $(v_0,\tht_0)\in D(A_1)\times D(A_2)$, there
 exists a unique solution $(v,\tht)$ of \eqref{e:v.n}-\eqref{e:ic.n} such that
\begin{equation*}
 (v,\tht)\in L^\infty(0, +\infty; D(A_1)\times D(A_2)).
\end{equation*}
 Moreover, there exists a bounded absorbing ball for the solutions of
 \eqref{e:v.n}-\eqref{e:ic.n} in $D(A_1)\times D(A_2)$.

\end{thrm}

\begin{prf}

 Uniqueness of such solutions follows immediately from the uniqueness of
 strong solutions.

 Notice that for any $(v_0,\tht_0)\in D(A_1)\times D(A_2)$, we have
 $(v_t(0), \tht_t(0))\in H$. Therefore, by Theorem~\ref{t:t}, we have
\[ (v_t,\tht_t)\in L^\infty(0, \infty; H)\cap L^2(0, \infty; V). \]
 Moreover, there exists a bounded absorbing ball for $(\prt_tv,\prt_t\tht)$ in
 space $H$.

 Now, we prove existence of a bounded absorbing set in $\R_+$ for both
 $\|v\|_{H^2}$ and $\|\tht\|_{H^2}$ and the uniform boundedness of $\|v\|_{H^2}$
 and $\|\tht\|_{H^2}$ for $t>0$.

{\bf Step 1}. Take inner product of \eqref{e:v.n} with $-\prt_z^2 v$ and use
 \eqref{e:bc.v} to obtain
\begin{equation}
\label{e:vz.h1}
\begin{split}
 \nu_1\|\nbl v_z\|_2^2 +& \nu_1\|v_{zz}\|_2^2\\
 &= \lra{v_t + (v\cdot\nbl)v + wv_z + \nbl p_s - \int_{-h}^z\nbl\tht
    + fv^\bot}{v_{zz}}.
\end{split}
\end{equation}
 First, by boundary condition \eqref{e:bc.v}, we have
\[ \lra{\nbl p_s}{v_{zz}} = \int_D \nbl p_s\int_{-h}^0 v_{zz}
 = \int_D \nbl p_s \lf(v_z\big|_{-h}^0\rt) = 0, \]
\[ -\lra{\int_{-h}^z\nbl\tht}{v_{zz}}
 =-\int_D\lf(\int_{-h}^z\nbl\tht\rt)\cdot v_z\big|_{z=-h}^0
  +\int_D\int_{-h}^0\nbl\tht\cdot v_z =\lra{\nbl\tht}{v_z},\]
 and, due to the fact that $f_z=0$,
\[ \lra{f v^\bot}{v_{zz}}=\int_D f \int_{-h}^0 v^\bot \cdot v_{zz}
 =\int_D f v^\bot v_z\big|_{z=-h}^0 - \int_\Omg f v_z^\bot\cdot v_z =0. \]
 Second,
\begin{equation*}
\begin{split}
 \lra{(v\cdot\nbl)v}{v_{zz}} =& \int_D [(v\cdot\nbl)v]\cdot v_z\big|_{z=-h}^0
 -\int_\Omg [(v_z\cdot\nbl)v+(v\cdot\nbl) v_z]\cdot v_z\\
=&-\int_\Omg [(v_z\cdot\nbl)v]\cdot v_z -\int_\Omg [(v\cdot\nbl) v_z]\cdot v_z\\
=&-\int_\Omg [(v_z\cdot\nbl)v]\cdot v_z+\hf \int_\Omg(\nbl\cdot v) |v_z|^2,
 \end{split}
\end{equation*}
 where in the last equality we have used the following computation:
\begin{equation*}
\begin{split}
 \int_\Omg [(v\cdot\nbl) v_z]\cdot v_z
=&\int_\Omg \sum_{i=1}^2 v_i\prt_i\lf(\frac{|v_z|^2}{2}\rt)\\
=&\int_{-h}^0 \hf \lf[\int_{\prt D}v\cdot n |v_z|^2
  - \int_{D}(\nbl\cdot v) |v_z|^2 \rt]\\
=& -\hf \int_\Omg(\nbl\cdot v) |v_z|^2.
\end{split}
\end{equation*}
 Third, since $w|_{z=-h}=w|_{z=0}=0$, we have
\begin{equation*}
\begin{split}
 \lra{wv_z}{v_{zz}} =& \hf \int_\Omg w \lf(|v_z|^2\rt)_z
 =\hf\int_D \lf(w |v_z|^2\big|_{z=-h}^{z=0} - \int_{-h}^0 w_z|v_z|^2 \rt)\\
 =& \hf \int_\Omg (\nbl\cdot v)|v_z|^2.
\end{split}
\end{equation*}
 Therefore,
\begin{equation*}
 \lra{(v\cdot\nbl)v}{v_{zz}} + \lra{wv_z}{v_{zz}}
= -\int_\Omg [(v_z\cdot\nbl)v]\cdot v_z + \int_\Omg(\nbl\cdot v) |v_z|^2
\end{equation*}
 Plugging the above computations into \eqref{e:vz.h1} and noticing that
 $v_z|_{z=0}=0$, we obtain
\begin{equation*}
\begin{split}
 \nu_1\|\nbl v_z\|_2^2 &+ \mu_1\|v_{zz}\|_2^2\\
&\ls \|v_t\|_2\|v_{zz}\|_2 + \|\nbl\tht\|_2\|v_z\|_2
 -\int_\Omg [(v_z\cdot\nbl)v]\cdot v_z + \int_\Omg(\nbl\cdot v) |v_z|^2\\
&\ls \|v_t\|_2\|v_{zz}\|_2 + \|\nbl\tht\|_2\|v_z\|_2
 + C\|\nbl v\|_2\|v_z\|_4^2\\
&\ls C_{\veps}\|v_t\|_2^2 +\veps\|v_{zz}\|_2^2 + \|\nbl\tht\|_2\|v_z\|_2
 + C\|\nbl v\|_2 \|v_z\|_2^\hf\|\nbl_3 v_z\|_2^{\frac{3}{2}}\\
&\ls C_{\veps}\lf(\|v_t\|_2^2 + \|\nbl v\|_2^4\|v_z\|_2^2\rt) + \|\nbl\tht\|_2\|v_z\|_2
 + 2\veps \|\nbl_3 v_z\|_2^2\\
&=C_{\veps}\lf(\|v_t\|_2^2 + \|\nbl v\|_2^4 \|v_z\|_2^2\rt) + \|\nbl\tht\|_2\|v_z\|_2
 + 2\veps (\|\nbl v_z\|_2^2 +\|v_{zz}\|_2^2).
\end{split}
\end{equation*}
 Thus, choosing $\veps>0$ sufficiently small yields
\begin{equation*}
 \|\nbl v_z\|_2^2 + \|v_{zz}\|_2^2\\
 \ls C (\|v_t\|_2^2 + \|\nbl v\|_2^4\|v_z\|_2^2 + \|\nbl\tht\|_2\|v_z\|_2).
\end{equation*}
 This gives the uniform boundedness of $\|v_z\|_{H^1}$ and a bounded absorbing
 set for $v_z$ in $H^1$.

 {\bf Step 2}. Take the inner product of \eqref{e:v.n} with $-\Dlt v$ and use
 \eqref{e:bc.v} to obtain
\begin{equation*}
\begin{split}
 \nu_1\|\Dlt v\|_2^2 + \mu_1\|\nbl v_z\|_2^2
 =& \lra{v_t}{\Dlt v} + \lra{(v\cdot\nbl)v}{\Dlt v} + \lra{wv_z}{\Dlt v}\\
   & -\lra{\int_{-h}^z\nbl\cdot\tht}{\Dlt v} +\lra{fv^\bot}{\Dlt v}.
\end{split}
\end{equation*}
 Therefore,
\begin{equation}
\label{v.nbl.h1}
\begin{split}
\|\Dlt v\|_2^2 + \|\nbl v_z\|_2^2
 \ls C(\|v_t\|_2^2 + \|(v\cdot\nbl)v\|_2^2 +\|wv_z\|_2^2
     +\|\nbl\tht\|_2^2 +\|v\|_2^2).
\end{split}
\end{equation}
 Notice that
\begin{equation*}
\begin{split}
 \|(v\cdot\nbl)v\|_2^2
\ls& \|v\|_6^2\|\nbl v\|_3^2\ls C\|v\|_{H^1}^2\|\nbl v\|_2\|\Dlt v\|_2\\
\ls& C_\veps\|v\|_{H^1}^6+ \veps\|\nbl v\|_{H^1}^2\\
=& C_\veps\|v\|_{H^1}^6+ \veps\|\nbl v\|_2^2
 + \veps(\|\Dlt v\|_2^2+\|\nbl v_z\|_2^2),
\end{split}
\end{equation*}
 and that by Lemma~\ref{l:ju.1}, we have
\begin{equation*}
\begin{split}
 \|wv_z\|_2^2
\ls& C\|\nbl v\|_2\|\nbl v\|_{H^1}\|v_z\|_2\|v_z\|_{H^1}\\
\ls& C_\veps\|\nbl v\|_2^2\|v_z\|_2^2\|v_z\|_{H^1}^2 +\veps \|\nbl v\|_{H^1}^2\\
=&C_\veps\|\nbl v\|_2^2\|v_z\|_2^2\|v_z\|_{H^1}^2 +\veps \|\nbl v\|_2^2
 +\veps (\|\Dlt v\|_2^2+\|\nbl v_z\|_2^2).
\end{split}
\end{equation*}
 Choosing $\veps>0$ sufficiently small, we derive from \eqref{v.nbl.h1} that
\begin{equation}
\label{e:v.dlt}
\begin{split}
 \|\Dlt v\|_2^2 + \|\nbl v_z\|_2^2
\ls& C(\|v_t\|_2^2 + \|v\|_{H^1}^6 + \|\nbl v\|_2^2\\
  &+ \|\nbl v\|_2^2\|v_z\|_2^2\|v_z\|_{H^1}^2 + \|\nbl\tht\|_2^2+ \|v\|_2^2).
\end{split}
\end{equation}
 Notice that there exists a bounded absorbing set in $\R_+$ for each term on
 the right-hand side of \eqref{e:v.dlt} and that each term on the right-hand
 side of \eqref{e:v.dlt} is uniformly bounded for $t>0$. Therefore, the same
 is true for $\|\Dlt v\|_2$, and we can conclude the same for $\|v\|_{H^2}$.

 {\bf Step 3}. Take the inner product of \eqref{e:t.n} with $-\prt_z^2\tht$ and
 use \eqref{e:bc.t} to obtain
\begin{equation}
\label{e:tz.h1}
 \nu_2\|\nbl\tht_z\|_2^2 + \alf\nu_2\|\nbl\tht|_{z=0}\|_2^2
 + \mu_2\|\tht_{zz}\|_2^2
 = \lra{\tht_t + v\cdot\nbl\tht+w\tht_z + Q}{\tht_{zz}}.
\end{equation}
 First,
\begin{equation*}
\begin{split}
 \lra{v\cdot\nbl\tht}{\tht_{zz}}
&\ls \|v\|_\infty\|\nbl\tht\|_2\|\tht_{zz}\|_2
 \ls C\|v\|_{H^2}\|\nbl\tht\|_2\|\tht_{zz}\|_2\\
&\ls C_\veps\|v\|_{H^2}^2\|\nbl\tht\|_2^2 + \veps\|\tht_{zz}\|_2^2.
\end{split}
\end{equation*}
 Second, using boundary condition $w|_{z=0}=w|_{z=-h}=0$, we have
\[  \lra{w\tht_z}{\tht_{zz}} = \hf\int_D\int_{-h}^0 w ((\tht_z)^2)_z
 = -\hf\int_\Omg w_z (\tht_z)^2=\hf\int_\Omg (\nbl\cdot v)(\tht_z)^2. \]
 Thus, noticing that $\tht_z|_{z=-h}=0$,
\begin{equation*}
\begin{split}
\lra{w\tht_z}{\tht_{zz}} &\ls \hf\|\nbl v\|_2 \|\tht_z\|_4^2
 \ls C\|\nbl v\|_2\|\tht_z\|_2^\hf\|\nbl_3\tht_z\|_2^{\frac{3}{2}}\\
 &\ls C_\veps\|\nbl v\|_2^4\|\tht_z\|_2^2 + \veps \|\nbl_3\tht_z\|_2^2.
\end{split}
\end{equation*}
 Plugging the above estimates into \eqref{e:tz.h1} and choosing sufficiently
 small $\veps>0$ then yields
\[ \|\nbl\tht_z\|_2^2 + \|\nbl\tht|_{z=0}\|_2^2
 + \|\tht_{zz}\|_2^2
  \ls C(\|\tht_t\|_2^2 + \|v\|_{H^2}^2\|\nbl\tht\|_2^2
  + \|\nbl v\|_2^4\|\tht_z\|_2^2 + \|Q\|_2^2). \]
 This finishes the proof of the uniform boundedness of $\|\tht_z\|_{H^1}$ and
 the existence of a bounded absorbing set of $\tht_z$ in $H^1$.

 {\bf Step 4}. Taking the inner product of \eqref{e:t.n} with $-\Dlt \tht$ and
 using \eqref{e:bc.t}, we obtain
\begin{equation*}
\begin{split}
 \nu_2\|\Dlt \tht\|_2^2 + \mu_2\|\nbl \tht_z\|_2^2
 =& \lra{\tht_t}{\Dlt\tht} + \lra{v\cdot\nbl\tht}{\Dlt\tht}
 + \lra{w\tht_z}{\Dlt\tht}\\
 \ls& C(\|\tht_t\|_2 + \|v\cdot\nbl\tht\|_2+\|w\tht_z\|_2)\|\Dlt\tht\|_2.
\end{split}
\end{equation*}
 Therefore, by Lemma~\ref{l:ju.1}, we obtain
\begin{equation*}
\begin{split}
 \|\Dlt \tht\|_2^2 &+ \|\nbl \tht_z\|_2^2\\
&\ls C(\|\tht_t\|_2^2 + \|v\cdot\nbl\tht\|_2^2+\|w\tht_z\|_2^2)\\
&\ls C(\|\tht_t\|_2^2 + \|v\|_{L^\infty}^2\|\nbl\tht\|_2^2
 + \|\nbl v\|_2\|\nbl v\|_{H^1}\|\tht_z\|_2\|\tht_z\|_{H^1})\\
&\ls C(\|\tht_t\|_2^2+ \|v\|_{H^2}^2\|\nbl\tht\|_2^2
 +\|\nbl v\|_2\|\nbl v\|_{H^1}\|\tht_z\|_2\|\tht_z\|_{H^1}),
\end{split}
\end{equation*}
 from which follows the existence of a bounded absorbing set for
 $\|\tht\|_{H^2}$ in $\R_+$ and the uniform boundedness of $\|\tht(t)\|_{H^2}$
 for $t>0$.

\end{prf}

\begin{rmrk}
 Following the ideas of \cite{j;t:2014}, Theorem~\ref{t:h2} can be used to give
 a different proof of the main result of \cite{j:pedim} on the finiteness of
 Hausdorff and fractal dimensions of the global attractor $\A$ for the strong
 solutions.
\end{rmrk}

\section{Global existence of $z$-weak Solutions}

\label{s:zw.e}

 Existence and uniqueness of $z$-weak solutions, global in time, was first
 proved in \cite{tachim:2010} for the case with boundary conditions
 \eqref{e:bc.v} and \eqref{e:bc.t}, under an extra assumption
 $(v_0, \tht_0)\in L^6$. As pointed out in \cite{tachim:2010} that this
 assumption can be replaced by a slightly weaker assumption that
 $(\tld{v_0},\tld{\tht_0})\in L^6$, where the notion of $\tld{\phi}$ is defined
 as the difference of a function $\phi$ on $\Omg$ and its vertical average:
\[ \tld{\phi}(x,y,z) := \phi(x,y,z) - \frac{1}{h}\int_{-h}^0\phi(x,y,z)\ dz. \]
 It seems to have been an open problem since the work of \cite{tachim:2010}
 that whether or not existence and uniqueness of $z$-weak solutions are valid
 for 3D Primitive equations, i.e. whether or not the extra assumption
 $(\tld{v_0},\tld{\tht_0})\in L^6$ can be eliminated or not. In this and next
 section, we resolve this problem with a positive confirmation. Our first
 result is the global existence of $z$-weak solutions as presented in the
 following Theorem~\ref{t:zw.e}. In next section, we deal with uniqueness of
 $z$-weak solutions.

 Now, we state our second main result of this article in the following:
\begin{thrm}
\label{t:zw.e}
 Suppose $Q\in L^2(\Omg)$. For any $(v_0,\tht_0),(\prt_zv_0,\prt_z\tht_0)\in H$,
 there exists a weak solution $(v,\tht)$ of \eqref{e:v.n}-\eqref{e:ic.n}
 such that
\begin{equation*}
 (v_z,\tht_z)\in L^\infty(0, +\infty; H)\cap L^2(0,\infty; V).
\end{equation*}
 Moreover, there exists a bounded absorbing ball for $(v_z,\tht_z)$ in $H$.

\end{thrm}

 \begin{prf}
 We prove the theorem in three steps:

{\bf Step 1}. Take inner product of \eqref{e:v.n} with $-\prt_z^2 v$ and use
 \eqref{e:bc.v} to obtain
\begin{equation}
\label{e:vz}
\begin{split}
 \hf\frac{d}{dt} \|v_z\|_2^2 +&\nu_1\|\nbl v_z\|_2^2+\nu_1\|v_{zz}\|_2^2\\
 &= \lra{(v\cdot\nbl)v +wv_z +\nbl p_s -\int_{-h}^z\nbl\tht+ fv^\bot}{v_{zz}},
\end{split}
\end{equation}
 where we have used the computation:
\[  -\int_\Omg v_t\cdot v_{zz} = -\int_D v_t\cdot v_z\Big|_{z=-h}^0  +
 \int_\omg v_{zt}\cdot v_z = \hf\frac{d}{dt} \|v_z\|_2^2.\]
 Moreover, by the few equalities in the previous subsection, and notice that
\begin{equation*}
\begin{split}
 -\int_\Omg \int_{-h}^z \nbl\tht\cdot v_{zz} = &
  -\int_D \lf(\int_{-h}^z\nbl\tht\rt)\cdot v_z \big|_{z=-h}^0
  +\int_\Omg \nbl\tht\cdot v_z \\
 =& \int_{-h}^0\int_{\prt D}\tht \frac{\prt v_z}{\prt n}
    -\int_\Omg \tht \nbl\cdot v_z =-\int_\Omg \tht \nbl\cdot v_z,
\end{split}
\end{equation*}
 we have
\begin{equation}
\label{e:vz.1}
\begin{split}
 \hf\frac{d}{dt} \|v_z\|_2^2 +&\nu_1\|\nbl v_z\|_2^2+\nu_1\|v_{zz}\|_2^2\\
 & = \int_\Omg \lf[(\nbl\cdot v)|v_z|^2-((v_z\cdot\nbl)v)\cdot v_z\rt]
 - \int_\Omg\tht\nbl\cdot v_z\\
 &= \int_\Omg \lf[(v_z\cdot \nbl^\bot) v^\bot\rt] \cdot v_z
 - \int_\Omg\tht\nbl\cdot v_z,
\end{split}
\end{equation}
 where
\[ v^\bot = (v_2, -v_1), \quad \nbl^\bot = (\prt_2, -\prt_1). \]
 By Lemma~\ref{l:ju.2}, we have
\begin{equation*}
\begin{split}
&\int_\Omg \lf[(v_z\cdot \nbl^\bot) v^\bot\rt] \cdot v_z\\
\ls& C \|\nbl v\|_2^\hf (\|\nbl v\|_2+\|\nbl v_z\|_2)^\hf
 \|v_z\|_2(\|v_z\|_2+\|\nbl v_z\|_2)\\
\ls& C(\|\nbl v\|_2 +\|\nbl v\|_2^\hf\|\nbl v_z\|_2^\hf)
 (\|v_z\|_2^2+\|v_z\|_2\|\nbl v_z\|_2) 
\end{split}
\end{equation*}
 By H\"older's inequality, we have
\begin{equation*}
\begin{split}
 &(\|\nbl v\|_2 +\|\nbl v\|_2^\hf\|\nbl v_z\|_2^\hf)
 (\|v_z\|_2^2+\|v_z\|_2\|\nbl v_z\|_2)\\
 =& \|\nbl v\|_2\|v_z\|_2^2 + \|\nbl v\|_2^\hf\|\nbl v_z\|_2^\hf\|v_z\|_2^2\\
  &+\|\nbl v\|_2\|v_z\|_2\|\nbl v_z\|_2
   +\|\nbl v\|_2^\hf\|v_z\|_2\|\nbl v_z\|_2^{\frac{3}{2}} \\
\ls& C_\veps\lf[(\|\nbl v\|_2^2+1)\|v_z\|_2^4 + 1\rt] + \veps\|\nbl v_z\|_2^2,
\end{split}
\end{equation*}
 and
\[ \lf|\int_\Omg\tht\nbl\cdot v_z\rt|
  \ls C_\veps\|\tht\|_2^2 + \veps\|\nbl v_z\|_2^2.\]
 Plugging the above estimates into \eqref{e:vz.1} and choosing  $\veps>0$
 sufficiently small yields
\begin{equation}
\label{e:vz.2}
 \frac{d}{dt}\|v_z\|_2^2 + \|\nbl v_z\|_2^2 + \|v_{zz}\|_2^2
 \ls  C\lf[ (\|\nbl v\|_2^2+1) \|v_z\|_2^4 + 1+ \|\tht\|_2^2\rt].
\end{equation}
 Therefore,
\[ \frac{d}{dt}\lf(\|v_z\|_2^2+1\rt)
   \ls C\lf(\|\nbl v\|_2^2+\|\tht\|_2^2+1\rt)(\|v_z\|_2^2+1)^2. \]
Let $z(t)=\|v_z\|_2^2 + 1$. Then, there exists a $t_1>0$, such that
\[ Cz(0)\int_0^{t_1}(\|\nbl v(\tau)\|_2^2+\|\tht(\tau)\|_2^2+1)d\tau < 1, \]
 and 
\[ z(t) \ls
 \frac{z(0)}{1 -z(0)\int_0^tC(\|\nbl v(\tau)\|_2^2+\|\tht(\tau)\|_2^2+1)d\tau},
 \quad \forall\ t\in[0,t_1]. \]
 This proves uniform boundedness of $\|v_z(t)\|_2 $ for $t\in[0, t_1]$. Then,
 by \eqref{e:vz.2}, we have uniform boundedness of
\[ \int_0^t\lf(\|\nbl v_z(\tau)\|_2^2+\|v_{zz}(\tau)\|_2^2\rt)\ d\tau, \quad
 \forall\ t\in[0, t_1]. \]

 {\bf Step 2}. Take inner product of \eqref{e:t.n} with $-\tht_{zz} $to obtain:
\begin{equation}
\label{e:tz}
\begin{split}
 \hf\frac{d}{dt}\lf(\|\tht_z\|_2^2 + \alf \|\tht|_{z=0}\|_2^2\rt)
 &+\nu_2\lf(\|\nbl\tht_z\|_2^2 + \alf\|\nbl\tht|_{z=0}\|_2^2\rt)
 + \mu_2\|\tht_{zz}\|_2^2\\
 =& \lra{Q}{\tht_{zz}} + \lra{v\cdot\nbl\tht+w\tht_z}{\tht_{zz}}.
\end{split}
\end{equation}
 In deriving \eqref{e:tz}, we have used boundary conditions \eqref{e:bc.v}
 and \eqref{e:bc.t} in the following calculations:
\[ -\int_\Omg \tht_t\tht_{zz} = -\int_D\tht_t\tht_z\big|_{z=-h}^0
  + \int_\Omg\tht_{zt}\tht_z
 = \hf\frac{d}{dt}\lf(\|\tht_z\|_2^2 + \alf\|\tht|_{z=0}\|_2^2\rt), \]
\begin{equation*}
\begin{split}
 \int_\Omg \Dlt\tht\cdot\tht_{zz} =& \int_D\Dlt\tht\cdot\tht_z\big|_{z=-h}^0 -
 \int_\Omg \Dlt\tht_z\cdot\tht_z
= -\alf\int_D\Dlt\tht\cdot\tht|_{z=0} -\int_\Omg \Dlt\tht_z\cdot\tht_z\\
=&-\alf\int_{\prt D} (n\cdot\nbl\tht)\tht|_{z=0}
 +\alf\int_{\prt D} \lf|\nbl\tht|_{z=0}\rt|^2\\
 &-\int_{-h}^0\int_{\prt D} (n\cdot\nbl\tht_z)\tht_z + \int_\Omg|\nbl\tht_z|^2\\
=& \alf\|\nbl\tht|_{z=0}\|_2^2 + \|\nbl\tht_z\|_2^2.
\end{split}
\end{equation*}
 Now we reformulate the right-hand side of equation \eqref{e:tz}.\\
 First of all, since $w|_{z=-h}=w|_{z=0}=0$, we have
\begin{equation*}
\begin{split}
\lra{w\tht_z}{\tht_{zz}} =& \hf \int_\Omg w (\tht_z^2)_z\\
 =& \hf \lf(\int_D w \tht_z^2\big|_{z=-h}^0 - \int_\Omg w_z \tht_z^2 \rt)
 = \hf \int_\Omg (\nbl\cdot v)\tht_z^2.
\end{split}
\end{equation*}
 Second,
\begin{equation*}
\begin{split}
 \lra{(v\cdot\nbl)\tht}{\tht_{zz}}
&= \int_D (v\cdot\nbl\tht)\tht_z\big|_{z=-h}^0
 - \int_\Omg\lf[(v_z\cdot\nbl\tht)\tht_z + (v\cdot\nbl\tht_z)\tht_z\rt] \\
&=-\alf\int_D (v\cdot\nbl\tht)\tht|_{z=0}
  -\int_\Omg (v_z\cdot\nbl\tht)\tht_z - \hf\int_\Omg v\cdot\nbl(\tht_z^2)\\
&=\frac{\alf}{2}\int_D \tht^2 \nbl v|_{z=0}
 -\int_\Omg (v_z\cdot\nbl\tht)\tht_z +\hf\int_\Omg(\nbl\cdot v)\tht_z^2.
\end{split}
\end{equation*}
 Thus,
\[  \lra{v\cdot\nbl\tht+w\tht_z}{\tht_{zz}}
 =\frac{\alf}{2}\int_D \tht^2 \nbl v|_{z=0}
 -\int_\Omg (v_z\cdot\nbl\tht)\tht_z + \int_\Omg(\nbl\cdot v)\tht_z^2.\]
 Therefore,
\begin{equation}
\label{e:tz.1}
\begin{split}
 \hf\frac{d}{dt}&\lf(\|\tht_z\|_2^2 + \alf\|\tht|_{z=0}\|_2^2\rt)
 +\nu_2\lf(\|\nbl\tht_z\|_2^2 + \alf\|\nbl\tht|_{z=0}\|_2^2\rt)
 + \mu_2\|\tht_{zz}\|_2^2\\
 =& \lra{Q}{\tht_{zz}} + \frac{\alf}{2}\int_D \tht^2 \nbl v|_{z=0}
 -\int_\Omg (v_z\cdot\nbl\tht)\tht_z +\int_\Omg(\nbl\cdot v)\tht_z^2\\
 =:& I_1 + I_2 + I_3 + I_4.
\end{split}
\end{equation}
 Now, we estimate $I_i$'s term by term. First, $I_1$ can be easily estimated
 as following:
\[ I_1 \ls C_\veps\|Q\|_2^2 + \veps\|\tht_{zz}\|_2^2. \]
 Next, since
\[ v|_{z=0}=v+\int_z^0v_z\ d\xi, \]
 we have, 
\begin{equation*}
\begin{split}
 h |v|_{z=0}|^2 &\ls 2\lf[ \int_{z=-h}^0 |v|^2
  + \int_{-h}^0 \lf(\int_z^0 |v_z|\ d\xi\rt)^2\ dz \rt]\\
&\ls 2\lf[ \int_{z=-h}^0 |v|^2 + h \int_{-h}^0|v_z|^2\ dz \rt]
\end{split}
\end{equation*}
 Thus,
\[ h \int_D |v|_{z=0}|^2 \ls 2 \lf[\int_\Omg |v|^2
  +  h \int_\Omg |v_z|^2\rt], \]
 that is
\begin{equation}
\label{e:v.z}
 \|v|_{z=0}\|_2^2 \ls \frac{2}{h} \|v\|_2^2 + 2\|v_z\|_2^2.
\end{equation}
 Therefore, by \eqref{e:v.z} and Gagliardo-Nirenberg-Sobolev inequality,
 we have
\begin{equation*}
\begin{split}
 I_2 \ls& \frac{\alf}{2} \|(\nbl\cdot v) |_{z=0}\|_2\|\tht |_{z=0}\|_4^2\\
\ls& C (\|\nbl\cdot v\|_2 +\|\nbl\cdot v_z\|_2) \|\tht|_{z=0}\|_2
   (\|\tht |_{z=0}\|_2 + \|\nbl\tht|_{z=0}\|_2)\\
\ls& C(\|\nbl v\|_2+\|\nbl v_z\|_2)\|\tht|_{z=0}\|_2^2\\
   & + C(\|\nbl v\|_2+\|\nbl v_z\|_2)\|\tht|_{z=0}\|_2\|\nbl\tht|_{z=0}\|_2\\
\ls& C(\|\nbl v\|_2+\|\nbl v_z\|_2)\|\tht|_{z=0}\|_2^2\\
   &+ C_\veps(\|\nbl v\|_2+\|\nbl v_z\|_2)^2\|\tht|_{z=0}\|_2^2
 +\veps\|\nbl\tht|_{z=0}\|_2^2\\
\ls& C_\veps(1 +\|\nbl v\|_2^2+\|\nbl v_z\|_2^2)\|\tht|_{z=0}\|_2^2
 +\veps\|\nbl\tht|_{z=0}\|_2^2.
\end{split}
\end{equation*}
 By Lemma~\ref{l:ju.2}, we have
\begin{equation*}
\begin{split}
 I_3 &\ls C\|v_z\|_2^\hf(\|v_z\|_2+\|\nbl v_z\|_2)^\hf
 \|\nbl\tht\|_2^\hf(\|\nbl\tht\|_2+\|\nbl\tht_z\|_2)^\hf\\
&\text{\hspace{.5in}}\times\|\tht_z\|_2^\hf(\|\tht_z\|_2+\|\nbl\tht_z\|_2)^\hf\\
&= C(\|v_z\|_2+\|v_z\|_2^\hf\|\nbl v_z\|_2^\hf)
 (\|\nbl\tht\|_2+\|\nbl\tht\|_2^\hf\|\nbl\tht_z\|_2^\hf)\\
&\text{\hspace{.5in}}\times(\|\tht_z\|_2+\|\tht_z\|_2^\hf\|\nbl\tht_z\|_2^\hf).
\end{split}
\end{equation*}
 Notice that
\begin{equation*}
\begin{split}
&(\|v_z\|_2+\|v_z\|_2^\hf\|\nbl v_z\|_2^\hf)
 (\|\nbl\tht\|_2+\|\nbl\tht\|_2^\hf\|\nbl\tht_z\|_2^\hf)
 (\|\tht_z\|_2+\|\tht_z\|_2^\hf\|\nbl\tht_z\|_2^\hf)\\
=&\|v_z\|_2\|\nbl\tht\|_2\|\tht_z\|_2 +
 \|v_z\|_2\|\nbl\tht\|_2\|\tht_z\|_2^\hf\|\nbl\tht_z\|_2^\hf\\
 &+\|v_z\|_2\|\nbl\tht\|_2^\hf\|\nbl\tht_z\|_2^\hf\|\tht_z\|_2
 +\|v_z\|_2\|\nbl\tht\|_2^\hf\|\nbl\tht_z\|_2^\hf
  \|\tht_z\|_2^\hf\|\nbl\tht_z\|_2^\hf\\
 &+\|v_z\|_2^\hf\|\nbl v_z\|_2^\hf\|\nbl\tht\|_2\|\tht_z\|_2
  +\|v_z\|_2^\hf\|\nbl v_z\|_2^\hf\|\nbl\tht\|_2
  \|\tht_z\|_2^\hf\|\nbl\tht_z\|_2^\hf\\
 &+\|v_z\|_2^\hf\|\nbl v_z\|_2^\hf\|\nbl\tht\|_2^\hf\|\nbl\tht_z\|_2^\hf
  \|\tht_z\|_2\\
 &+\|v_z\|_2^\hf\|\nbl v_z\|_2^\hf\|\nbl\tht\|_2^\hf\|\nbl\tht_z\|_2^\hf
 \|\tht_z\|_2^\hf\|\nbl\tht_z\|_2^\hf\\
&=:J_1+\dots+J_8.
\end{split}
\end{equation*}
 and, by H\"older's inequality, we have
\begin{equation*}
\begin{split}
 J_1 &\ls \hf\|v_z\|_2^2\|\tht_z\|_2^2 + \hf\|\nbl\tht\|_2^2,\\
 J_2 &\ls C_\veps(\|v_z\|_2^4\|\tht_z\|_2^2 + \|\nbl\tht\|_2^2)
 + \veps\|\nbl\tht_z\|_2^2,\\
 J_3 &\ls C_\veps(\|v_z\|_2^2\|\tht_z\|_2^2 + \|\nbl\tht\|_2^2)
 + \veps\|\nbl\tht_z\|_2^2, \\
 J_4 &= \|v_z\|_2\|\nbl\tht\|_2^\hf
  \|\tht_z\|_2^\hf\|\nbl\tht_z\|_2\\
 &\ls C_\veps(\|v_z\|_2^4\|\tht_z\|_2^2 + \|\nbl\tht\|_2^2)
 + \veps\|\nbl\tht_z\|_2^2,\\
 J_5 &\ls \frac{1}{4}\lf(\|v_z\|_2^2+\|\nbl v_z\|_2^2\rt)\|\tht_z\|_2^2
           + \hf\|\nbl\tht\|_2^2,\\
 J_6 &\ls C_\veps(\|v_z\|_2^2\|\nbl v_z\|_2^2\|\tht_z\|_2^2 + \|\nbl\tht\|_2^2)
 + \veps\|\nbl\tht_z\|_2^2,\\
 J_7 &\ls C_\veps\lf[(\|v_z\|_2^2+\|\nbl v_z\|_2^2)\|\tht_z\|_2^2
 + \|\nbl\tht\|_2^2\rt] + \veps\|\nbl\tht_z\|_2^2,\\
 J_8 &=\|v_z\|_2^\hf\|\nbl v_z\|_2^\hf\|\nbl\tht\|_2^\hf
 \|\tht_z\|_2^\hf\|\nbl\tht_z\|_2\\
 &\ls C_\veps(\|v_z\|_2^2\|\nbl v_z\|_2^2\|\tht_z\|_2^2 + \|\nbl\tht\|_2^2)
 + \veps\|\nbl\tht_z\|_2^2.
\end{split}
\end{equation*}
 Thus, for $t\in[0, t_1]$,
\[ I_3 \ls C_\veps(1+\|\nbl v_z\|_2^2)\|\tht_z\|_2^2 +C_\veps\|\nbl\tht\|_2^2
 + C\veps\|\nbl\tht_z\|_2^2 .\]
 By Lemma~\ref{l:ju.2} again, we have
\begin{equation*}
\begin{split}
 I_4 \ls& C\|\nbl v\|_2^\hf(\|\nbl v\|_2+\|\nbl v_z\|_2)^\hf
 \|\tht_z\|_2(\|\tht_z\|_2 +\|\nbl\tht_z\|_2)\\
 \ls& C(\|\nbl v\|_2+\|\nbl v\|_2^\hf\|\nbl v_z\|_2^\hf)
       (\|\tht_z\|_2^2 +\|\tht_z\|_2\|\nbl\tht_z\|_2).
\end{split}
\end{equation*}
 By H\"older's inequality, we have
\begin{equation*}
\begin{split}
&(\|\nbl v\|_2+\|\nbl v\|_2^\hf\|\nbl v_z\|_2^\hf)
       (\|\tht_z\|_2^2 +\|\tht_z\|_2\|\nbl\tht_z\|_2)\\
=&\|\nbl v\|_2\|\tht_z\|_2^2 + \|\nbl v\|_2^\hf\|\nbl v_z\|_2^\hf\|\tht_z\|_2^2\\
 &+ \|\nbl v\|_2\|\tht_z\|_2\|\nbl\tht_z\|_2
 +\|\nbl v\|_2^\hf\|\nbl v_z\|_2^\hf\|\tht_z\|_2\|\nbl\tht_z\|_2\\
\ls& C_\veps(1 + \|\nbl v\|_2^2 + \|\nbl v_z\|_2^2)\|\tht_z\|_2^2
  +\veps \|\nbl\tht_z\|_2^2.
\end{split}
\end{equation*}
 Thus,
\[ I_4 \ls C_\veps(1 + \|\nbl v\|_2^2 + \|\nbl v_z\|_2^2)\|\tht_z\|_2^2
  +C\veps \|\nbl\tht_z\|_2^2.\]
 Therefore, for $t\in[0, t_1]$, we have
\begin{equation*}
\begin{split}
 \hf\frac{d}{dt}\lf(\|\tht_z\|_2^2 + \alf\|\tht|_{z=0}\|_2^2\rt)
&+ \nu_2\lf(\|\nbl\tht_z\|_2^2 + \alf\|\nbl\tht|_{z=0}\|_2^2\rt)
 + \mu_2\|\tht_{zz}\|_2^2\\
\ls & C_\veps\|Q\|_2^2
 + C_\veps (1+ \|\nbl v\|_2^2+\|\nbl v_z\|_2^2)\|\tht|_{z=0}\|_2^2\\
  &+ C_\veps(1+\|\nbl v\|_2^2+\|\nbl v_z\|_2^2)\|\tht_z\|_2^2
   + C_\veps \|\nbl\tht\|_2^2 \\
 &+ \veps \|\tht_{zz}\|_2^2 + \veps\|\nbl\tht|_{z=0}\|_2^2
 + C\veps\|\nbl\tht_z\|_2^2.
\end{split}
\end{equation*}
 Therefore, choosing sufficiently small $\veps>0$, we have, for $t\in[0, t_1]$,
\begin{equation}
\label{e:tz.2}
\begin{split}
&\frac{d}{dt}\lf(\|\tht_z\|_2^2 + \alf\|\tht|_{z=0}\|_2^2\rt)
 + \nu_2\|\nbl\tht_z\|_2^2 + \nu\alf\|\nbl\tht|_{z=0}\|_2^2
 + \mu_2\|\tht_{zz}\|_2^2\\
\ls& C(\|Q\|_2^2+\|\nbl\tht\|_2^2)
 + C(1 + \|\nbl v\|_2^2+\|\nbl v_z\|_2^2)
 (\|\tht_z\|_2^2+\alf\|\tht|_{z=0}\|_2^2).
\end{split}
\end{equation}
 Noticing that
\[ \int_0^{t_1}(1 + \|\nbl v\|_2^2+\|\nbl v_z\|_2^2 + \|\nbl\tht\|_2^2) dt
 < +\infty, \]
 and \eqref{e:v.z}, i.e.
\[ \lf\|\tht|_{z=0}\rt\|_2^2 \ls \frac{2}{h}\|\tht\|_2^2 + 2\|\tht_z\|_2^2, \]
 by Gronwall's inequality, we obtain from \eqref{e:tz.2} uniform boundedness of
\[ \|\tht_z(t)\|_2^2 + \alf\|\tht(t)|_{z=0}\|_2^2 \]
 for $t\in [0,t_1]$.

{\bf Step 3}. Now, we can use global existence, uniqueness and uniform
 boundedness of the strong solutions to obtain global existence and uniform
 boundedness of $(v_z,\tht_z)\in L^2$. Indeed, since $(v,\tht)$ is a weak
 solution, there exists a $t_0\in(0, t_1]$ such that
\[ \|v(t_0)\|_{V_1}<+\infty,\quad \|\tht(t_0)\|_{V_2} <+\infty. \]
 Therefore, both $\|v(t)\|_{V_1}$ and $\|\tht\|_{V_2}$ are uniformly bounded
 for $t\in[t_0,\infty)$. Moreover, there exists a bounded absorbing set for 
 $(v,\tht)$ in $V$. Therefore, both $\|v_z(t)\|_2$ and $\|\tht_z(t)\|_2$ are
 uniformly bounded for $t\gs 0$ and there is a bounded absorbing set for
 $(v_z,\tht_z)$ in $H$.
\end{prf}

\section{Uniqueness of $z$-weak Solutions}
\label{s:zw.u}

 Now we prove uniqueness of $z$-weak solutions. Indeed, we will prove the
 following Lipschitz continuity of the solutions in $H$ space.

\begin{thrm}
\label{t:zw.u}
 Suppose $Q\in L^2(\Omg)$, $(v_0,\tht_0),(\prt_zv_0,\prt_z\tht_0)\in H $, then
 exists a unique weak solution $(v,\tht)$ of \eqref{e:v.n}-\eqref{e:ic.n}, which
 depends continuously on the initial data in $H$.
 More generally, let $(v^{(i)}, \tht_i)$ be a weak solution of the system of
 Primitive Equations with initial data $(v^{(i)}(0), \tht_i(0))\in H$ with
 $i=1,2$. Suppose $\|\prt_zv^{(1)}(0)\|_2$ and $\|\prt_z\tht_1(0)\|_2$ are
 finite.
 Then, for any $T>0$, there exists a constant $C>0$ which may depends on
 $(v^{(1)}(0), \tht_1(0))$ and $T$, but not on $(v^{(2)}, \tht_2)$, such that,
 for all $t\in[0,T]$,
\begin{equation*}
\begin{split}
 \|(v^{(1)}(t)-v^{(2)}(t),\ & \tht_1(t)-\tht_2(t))\|_H\\
 &\ls C \|(v^{(1)}(0)-v^{(2)}(0),\ \tht_1(0)-\tht_2(0))\|_H.
\end{split}
\end{equation*}
\end{thrm}
\begin{prf}

 Denote
 $p_{s,i}$ as the corresponding pressure term in the equation of $v^{(i)}$,
\[ w_i = -\int_{-h}^z \nbl\cdot v^{(i)}(x,y,\xi)\ d\xi. \]
 and
\[ u=v^{(1)}-v^{(2)},\quad  \omg=w_1-w_2, \quad \zt=\tht_1-\tht_2,\quad
   q_s=p_{1,s}-p_{2,s}.\]
 Then
\begin{equation}
\label{e:u}
\begin{split}
 u_t -\nu_1\Dlt u -& \mu_1 u_{zz} + (v^{(2)}\cdot\nbl) u  + w_2 u_z\\
  +& (u\cdot\nbl)v^{(1)} + \omg v^{(1)}_z 
  + \nbl q_s - \int_{-h}^z\nbl\zt + fu^\bot =0.
\end{split}
\end{equation}
\begin{equation}
\label{e:z}
\begin{split}
 \zt_t-\nu_2\Dlt\zt-\mu_2\zt_{zz} + v^{(2)}\cdot\nbl\zt + w_2\zt_z
  + u\cdot\nbl\tht_1 + \omg \tht_{1,z} = 0.
\end{split}
\end{equation}
 Take inner product of \eqref{e:u} with $u$ to obtain
\begin{equation}
\label{e:u.2}
\begin{split}
 \hf \dt{\|u\|_2^2} +& \nu_1\|\nbl u\|_2^2 + \mu_1\|u_z\|_2^2\\
 &= - \int_\Omg \lf[(u\cdot\nbl)v^{(1)} + \omg v^{(1)}_z -\int_{-h}^z\nbl\zt
 \rt]\cdot u,
\end{split}
\end{equation}
 where we have used the following fact:
\[ \int_\Omg \lf[ (v^{(2)}\cdot\nbl) u + w_2 u_z\rt] \cdot u
  =\int_\Omg \nbl q_s \cdot u = \int_\Omg f u^\bot\cdot u = 0. \]
 Take inner product of \eqref{e:z} with $\zt$ to obtain
\begin{equation}
\label{e:z.2}
\begin{split}
 \hf \dt{\|\zt\|_2^2} +& \nu_2\|\nbl\zt\|_2^2 + \mu_2\|\zt_z\|_2^2
 + \alf\mu_2\|\zt|_{z=0}\|_2^2\\
 &= -\int_\Omg \lf[ (v\cdot\nbl)\tht_1 + \omg\tht_{1,z}\rt]\zt,
\end{split}
\end{equation}
 where we have used the following fact:
\[ \int_\Omg \lf[(v^{(2)}\cdot\nbl)\zt + w_2\zt_z\rt]\zt =0. \]
 First, it is easy to have that
\begin{equation*}
\begin{split}
 \int_\Omg \int_{-h}^z(\nbl\zt)\cdot u = -\int_\Omg\int_{-h}^z\zt(\nbl\cdot u)
 \ls C_\veps h^2\|\zt\|_2^2 + \veps\|\nbl u\|_2^2.
\end{split}
\end{equation*}
 Next, we use Lemma~\ref{l:ju.2} and H\"older's inequality to get the
 following estimate:
\begin{equation*}
\begin{split}
 \int_\Omg ((u\cdot\nbl) v^{(1)})\cdot u
 \ls& C \|u\|_2\|\nbl v^{(1)}\|_2^\hf(\|u\|_2+\|\nbl u\|_2)
  (\|\nbl v^{(1)}\|_2 + \|\nbl v^{(1)}_z\|_2)^\hf\\
 \ls& C_\veps\|\nbl v^{(1)}\|_2(\|\nbl v^{(1)}\|_2+\|\nbl v^{(1)}_z\|_2)
 \|u\|_2^2\\ &+\frac{\veps}{2}(\| u\|_2+\|\nbl u\|_2)^2\\
 \ls &C_\veps(\|\nbl v^{(1)}\|_2^2+\|\nbl v^{(1)}_z\|_2^2)\|u\|_2^2
 + \veps(\| u\|_2^2+\|\nbl u\|_2^2)\\
 \ls &C_\veps(1 + \|\nbl v^{(1)}\|_2^2+\|\nbl v^{(1)}_z\|_2^2)\|u\|_2^2
 +\veps\|\nbl u\|_2^2.
\end{split}
\end{equation*}
 Noticing that
\[ \|\omg\|_2 = \lf\|-\int_{-h}^z \nbl\cdot u\rt\|_2 \ls  h\|\nbl\cdot u\|_2,
  \quad \|\omg_z\|_2 = \|\nbl\cdot u\|_2, \]
 we use Lemma~\ref{l:ju.2} and H\"older's inequality again to obtain:
\begin{equation*}
\begin{split}
 \int_\Omg \omg v^{(1)}_z\cdot u
 \ls& C \|\omg\|_2^\hf\|v^{(1)}_z\|_2^\hf\|u\|_2^\hf
  (\|\omg\|_2+\|\omg_z\|_2)^\hf\\
  &\times  (\|v^{(1)}_z\|_2+\|\nbl v^{(1)}_z\|_2)^\hf
  (\|u\|_2+ \|\nbl u\|_2)^\hf\\
 =& C \|v^{(1)}_z\|_2^\hf (\|v^{(1)}_z\|_2+\|\nbl v^{(1)}_z\|_2)^\hf
   \|u\|_2^\hf \|\nbl u\|_2 (\|u\|_2+\|\nbl u\|_2)^\hf\\
 \ls& C_\veps(1 +\|v^{(1)}_z\|_2^4 +\|v^{(1)}_z\|_2^2\|\nbl v^{(1)}_z\|_2^2)
 \|u\|_2^2 + \veps\|\nbl u\|_2^2\\
 \ls& C_\veps(1 +\|\nbl v^{(1)}_z\|_2^2)\|u\|_2^2 + \veps\|\nbl u\|_2^2,
\end{split}
\end{equation*}
 where we have used uniform boundedness of $\|v^{(1)}_z\|_2$.

 Similar to the treatment of the above two nonlinear terms, we also have
\begin{equation*}
\begin{split}
\int_\Omg (u\cdot\nbl\tht_1) \zt
\ls& C\|u\|_2^\hf\|\nbl\tht_1\|_2^\hf\|\zt\|_2^\hf
      (\|u\|_2+\|\nbl u\|_2)^\hf\\
 &\times (\|\nbl\tht_1\|_2+\|\nbl\tht_{1,z}\|_2)^\hf
      ( \|\zt\|_2+\|\nbl\zt\|_2)^\hf\\
\ls& C_\veps\|\nbl\tht_1\|_2^2\|u\|_2^2 +
 C_\veps(\|\nbl\tht_1\|_2^2+\|\nbl\tht_{1,z}\|_2^2)\|\zt\|_2^2\\
 & + \veps(\|u\|_2^2+ \|\nbl u\|_2^2  + \|\zt\|_2^2+ \|\nbl\zt\|_2^2)\\
\ls& C_\veps(1+\|\nbl\tht_1\|_2^2+\|\nbl\tht_{1,z}\|_2^2)
 (\|u\|_2^2+\|\zt\|_2^2)\\
 &+ \veps(\|\nbl u\|_2^2+\|\nbl\zt\|_2^2),
\end{split}
\end{equation*}
 and
\begin{equation*}
\begin{split}
 \int_\Omg \omg\tht_{1,z}\zt
\ls& C\|\omg\|_2^\hf \|\tht_{1,z}\|_2^\hf\|\zt\|_2^\hf
      (\|\omg\|_2+\|\omg_z\|_2)^\hf\\ &\times
      (\|\tht_{1,z}\|_2+\|\nbl\tht_{1,z}\|_2)^\hf
      (\|\zt\|_2+\|\nbl\zt\|_2)^\hf\\
\ls& C\|\tht_{1,z}\|_2^\hf(\|\tht_{1,z}\|_2+\|\nbl\tht_{1,z}\|_2)^\hf
      \|\zt\|_2^\hf \|\nbl u\|_2(\|\zt\|_2+\|\nbl\zt\|_2)^\hf\\
\ls& C_\veps\|\tht_{1,z}\|_2^2(\|\tht_{1,z}\|_2^2+\|\nbl\tht_{1,z}\|_2^2)
  \|\zt\|_2^2 +\veps(\|\nbl u\|_2^2+\|\zt\|_2^2+\|\nbl\zt\|_2^2)\\
\ls& C_\veps(1+\|\tht_{1,z}\|_2^4+\|\tht_{1,z}\|_2^2\|\nbl\tht_{1,z}\|_2^2)
 \|\zt\|_2^2  +\veps(\|\nbl u\|_2^2 +\|\nbl\zt\|_2^2)\\
\ls& C_\veps(1+\|\nbl\tht_{1,z}\|_2^2)\|\zt\|_2^2 
  +\veps(\|\nbl u\|_2^2 +\|\nbl\zt\|_2^2).
\end{split}
\end{equation*}
where we have used uniform boundedness of $\|\tht_{1,z}\|_2$.

 Plugging the above estimates into \eqref{e:u.2} and \eqref{e:z.2} and
 adding up the two inequalities yields
\begin{equation*}
\begin{split}
 \hf \dt{(\|u\|_2^2+\|\zt\|_2^2)} +& \nu_1\|\nbl u\|_2^2 + \mu_1\|u_z\|_2^2
 + \nu_2\|\nbl\zt\|_2^2 + \mu_2(\|\zt_z\|_2^2+ \alf\|\zt|_{z=0}\|_2^2)\\
\ls& C_\veps(1 + \|\nbl v^{(1)}\|_2^2+\|\nbl v^{(1)}_z\|_2^2+
 \|\nbl\tht_1\|_2^2+\|\nbl\tht_{1,z}\|_2^2)\\
 &\times (\|u\|_2^2+\|\zt\|_2^2) + 4\veps(\|\nbl u\|_2^2+\|\nbl\zt\|_2^2).
\end{split}
\end{equation*}
 Choose $\veps>0$ sufficiently small, we have
\begin{equation}
\label{e:u.z}
\begin{split}
 \dt{(\|u\|_2^2+\|\zt\|_2^2)} +& \nu_1\|\nbl u\|_2^2 + \mu_1\|u_z\|_2^2
 + \nu_2\|\nbl\zt\|_2^2 + \mu_2(\|\zt_z\|_2^2 + \alf\|\zt|_{z=0}\|_2^2)\\
\ls& C_\veps(1 + \|\nbl v^{(1)}\|_2^2+\|\nbl v^{(1)}_z\|_2^2+
 \|\nbl\tht_1\|_2^2+\|\nbl\tht_{1,z}\|_2^2)\\
 &\times (\|u\|_2^2+\|\zt\|_2^2).
\end{split}
\end{equation}
 Noticing that, for a $z$-weak solution $(v^{(1)},\tht_1)$,
\[ \int_0^\infty (\|\nbl v^{(1)}_z(t)\|_2^2+\|\nbl\tht_{1,z}(t)\|_2^2)\ dt
  <+\infty, \]
 Integrating \eqref{e:u.z} with respect to $t$, after dropping the dissipation
 and diffusion terms, yields Lipschitz continuity of $(v,\tht)$ with respect to
 $(v_0,\tht_0)$ in space $H$, thus the uniqueness of $z$-weak solutions.

\end{prf}

{\bf Acknowledgment:}\\
The author wishes to thank Professor Roger Temam for consistent encouragement
and precious support of this research project.

\end{document}